\documentclass[reviewcopy]{elsart}
\usepackage{natbib}
\usepackage{graphicx}
\usepackage{amssymb}
\usepackage{amsmath}
\bibpunct{[}{]}{,}{n}{,}{,}	% see natbib.sty
\newcommand\assign{\mathbin{:=}}

\newcommand\vdp{\boldsymbol{\cdot}}
\begin{document}
\begin{frontmatter}
\date{2005 June 2}
\title{Exact calculation of Fourier series in
nonconforming spectral-element methods}%\thanksref{label1}}
\author{Aim\'e Fournier}%\corauthref{cor1}\thanksref{label2}}
\ead{fournier@ucar.edu}
\ead[url]{http://www.asp.ucar.edu/gtp/fournier}
\address{NCAR Institute for Mathematics Applied to Geosciences\\
Boulder CO 80307-3000 USA}
% \begin{abstract}
% Text of abstract
% \end{abstract}
\end{frontmatter}

% main text
\section{Usefulness of calculating Fourier series in the SEM}

In this note is presented a method, given nodal values
on multidimensional nonconforming spectral elements,
for calculating global Fourier-series coefficients.
%-----------------------------------------------------------------------
This method is ``exact'' in that given the approximation
inherent in the spectral-element method (SEM), no further approximation
is introduced that exceeds computer round-off error.
%-----------------------------------------------------------------------
The method is very useful when the SEM has yielded
an adaptive-mesh representation of a spatial function whose global
Fourier spectrum must be examined, e.g., in dynamically adaptive
fluid-dynamics simulations such as \citep{RFFP2005}.

\section{Derivation of an exact transform}
\label{mthdlg}

Suppose we have some functional problem in a spatial
domain $\mathbb{D}\assign[-\pi,\pi]^d$ (possibly including
toroidal geometry) and use coordinate transformations
\begin{equation}%-------------------------------------------------------
\vec\vartheta_k\quad\text{from}\quad
\vec\xi\in\mathbb{E}_0:=[-1,1]^d\quad\text{to}\quad\vec{x}\in\mathbb{E}_k
\label{lct}%------------------------------------------------------------
\end{equation}%---------------------------------------------------------
to partition $\mathbb{D}=\bigcup_{k=1}^K\mathbb{E}_k$ by
$K$ elements $\mathbb{E}_k\assign\vec\vartheta_k(\mathbb{E}_0)$ with
disjoint\footnote{$\overset{\bullet}{\mathbb{E}}_k\bigcap\overset{\bullet}{\mathbb{E}}_{k'}=\varnothing$
if $k\neq{k'}$} interiors.
%-----------------------------------------------------------------------
Typically the SEM approximates the exact solution by
its piecewise polynomial representation of degree $P$:
%-----------------------------------------------------------------------
\begin{equation}%-------------------------------------------------------
u^{\rm{ex}}(\vec{x})\approx
u(\vec{x})=\sum_{k=1}^K\sum_{\vec\jmath\in\mathbb{J}}
u_{\vec\jmath,k}\phi_{\vec\jmath,k}(\vec{x}),
\label{ppr}%------------------------------------------------------------
\end{equation}%---------------------------------------------------------
where $\mathbb{J}:=\{0,\ldots P\}^d$ indexes the values
$u_{\vec\jmath,k}\assign u(\vec{x}_{\vec\jmath,k})$ and nodes
$\vec{x}_{\vec\jmath,k}\assign\vec\vartheta_k(\vec\xi_{\vec\jmath})$
mapped from the $d$-dimensional Gauss-Lobatto-Legendre (GLL) quadrature
nodes $\xi^\alpha_{\vec\jmath}\assign\xi_{\jmath^\alpha}\in[-1,1]$,
\begin{equation}%-------------------------------------------------------
\phi_{\vec\jmath,k}(\vec{x})\assign\begin{cases}
\phi_{\vec\jmath}\circ\vec\vartheta^{-1}_k(\vec{x}),&\vec{x}\in\mathbb{E}_k\\
0,&\vec{x}\not\in\mathbb{E}_k\end{cases}
\label{mappedip}%-------------------------------------------------------
\end{equation}%---------------------------------------------------------
is the $\vec{x}_{\vec\jmath,k}$-interpolating piecewise-polynomial,
\begin{equation}%-------------------------------------------------------
\phi_{\vec\jmath}(\vec\xi)\assign
\prod_{\alpha=1}^d\phi_{\jmath^\alpha}(\xi^\alpha) \quad\text{and}\quad
\phi_j(\xi)\assign\sum_{p=0}^P\check{\phi}_{j,p}{\rm{L}}_{p}(\xi)
\label{GLLip}
\end{equation}%---------------------------------------------------------
are $\vec\xi_{\vec\jmath}$\:- and $\xi_j$-interpolating
polynomials, $\check{\phi}_{j,p}\equiv
w_j{\rm{L}}_{p}(\xi_j)/\sum_{j'=0}^Pw_{j'}{\rm{L}}_{p}(\xi_{j'})^2$
is a Legendre coefficient \citep[e.g.,][(B.3.15)]{DFM2002},
$\sqrt{p+\half}{\rm{L}}_p(\xi)$ is the orthonormal Legendre polynomial
of degree $p$ on $[-1,1]$ and $w_j$ is the GLL quadrature weight.
%-----------------------------------------------------------------------
In many cases a physically interesting quantity is the global
Fourier-series coefficient $\hat{u}_{\vec{q}}$ at integer
wavenumber components $q^\alpha$, usually approximated
by $M^d$-point trigonometric $d$-cubature in such manner as
\begin{align}%----------------------------------------------------------
\hat{u}_{\vec{q}}&\assign\frac{1}{(2\pi)^d}\int_{\mathbb{D}}
u(\vec{x})\e^{-{\rm{i}}\vec{q}\vdp\vec{x}}\d v(\vec{x})
\equiv\sum_{k=1}^K\sum_{\vec\jmath\in\mathbb{J}}
\hat{\phi}_{\vec\jmath,k,\vec{q}}u_{\vec\jmath,k},
\label{Fsc}\\%----------------------------------------------------------
\text{where}\quad\hat{\phi}_{\vec\jmath,k,\vec{q}}&=
\frac{1}{M^d}\sum_{\vec{m}\in\mathbb{M}}
\phi_{\vec\jmath,k}(\vec{x}_{\vec{m}})\e^{-{\rm{i}}\vec{q}\vdp\vec{x}_{\vec{m}}}
-\mathcal{E}_{\vec{q}}\phi_{\vec\jmath,k},
\label{aFsc}%-----------------------------------------------------------
\end{align}%------------------------------------------------------------
$\d v(\vec{x}):=\prod_{\alpha=1}^d\d x^\alpha$ is the volume
differential and $\mathbb{M}:=\{1,\ldots M\}^d$ indexes
trigonometric nodes $x^\alpha_{\vec{m}}\assign(2m^\alpha/M-1)\pi$.
%-----------------------------------------------------------------------
Note whenever $\mathbb{D}$ is adaptively repartitioned there
is an additional computation cost of $\mathcal{O}(M^d)$
per node to use \eqref{ppr} to provide in \eqref{aFsc} the
values $\phi_{\vec\jmath,k}(\vec{x}_{\vec{m}})$, as well as a
$d$-cubature error \citep[generalizing][theorem 4.7]{Boyd1989}
\begin{align*}%---------------------------------------------------------
\mathcal{E}_{\vec{q}}u\equiv\sum_{\vec{r}\in
\Zset^d\setminus\{\vec0\}}\hat{u}_{\vec{q}+M\vec{r}}
%\label{eFsc}-----------------------------------------------------------
\end{align*}%-----------------------------------------------------------
that in general converges no faster than
$\mathcal{O}(M^{-2})$, because $\mathbb{C}^1$
discontinuities of \eqref{ppr} across element boundaries
%$\overset{\circ}{\mathbb{E}}_k:=\mathbb{E}_k\setminus\overset{\bullet}{\mathbb{E}}_k$
cause $|\hat{u}_{\vec{q}}|$
to decay only as $\mathcal{O}(|\vec{q}|^{-2})$.
%-----------------------------------------------------------------------
We discover a more accurate method by
substituting \eqref{mappedip} into \eqref{Fsc} to yield
\begin{align*}%---------------------------------------------------------
\hat{\phi}_{\vec\jmath,k,\vec{q}}&=\frac{1}{(2\pi)^d}
\int_{\mathbb{E}_k}\e^{-{\rm{i}}\vec{q}\vdp\vec{x}}
\phi_{\vec\jmath}\circ\vec\vartheta^{-1}_k(\vec{x})\d v(\vec{x})\\
&\overset{\eqref{lct}}{=}\frac{1}{(2\pi)^d}\int_{\mathbb{E}_0}
\e^{-{\rm{i}}\vec{q}\vdp\vec\vartheta_k(\vec\xi)}
\phi_{\vec\jmath}(\vec\xi)
\left|\frac{\partial\vec\vartheta_k}{\partial\vec\xi}\right|\d
v(\vec\xi)\\%-----------------------------------------------------------
&\overset{\eqref{GLLip}}{=} \frac{1}{(2\pi)^d}\int_{\mathbb{E}_0}
\e^{-{\rm{i}}\vec{q}\vdp\vec\vartheta_k(\vec\xi)}
\left(\prod_{\alpha=1}^d\sum_{p=0}^P
\check{\phi}_{\jmath^\alpha,p}{\rm{L}}_p(\xi^\alpha)\right)
\left|\frac{\partial\vec\vartheta_k}{\partial\vec\xi}\right|\d
v(\vec\xi).
\end{align*}%-----------------------------------------------------------
In many applications, especially when $u$-structure rather
than domain geometry is guiding the mesh adaption, each
$\mathbb{E}_k$ is a $d$-parallelepiped with center $\vec{a}_k$
and $d$ legs $2\vec{h}^\alpha_k$, so we have an affinity
$\vec\vartheta_k(\vec\xi):=\vec{a}_k+\vec{\vec{h}}_k\vdp\vec\xi$,
where $\vec{h}^\alpha_k$ make up the columns of $\vec{\vec{h}}_k$.
%-----------------------------------------------------------------------
Then we obtain
$$%---------------------------------------------------------------------
\hat{\phi}_{\vec\jmath,k,\vec{q}}=
\frac{1}{(2\pi)^d}\left|\vec{\vec{h}}_k\right|
\e^{-{\rm{i}}\vec{q}\vdp\vec{a}_k} \prod_{\alpha=1}^d\sum_{p=0}^P
\check{\phi}_{\jmath^\alpha,p}\int_{-1}^1
\e^{-{\rm{i}}\vec{q}\vdp\vec{h}^\alpha_k\xi} {\rm{L}}_p(\xi)\d\xi.
$$%--------------------------------------------------------------------
Finally, recalling the classical identity \citep[e.g.,][exercise 12.4.9]{Arf85}
for the spherical Bessel function
${\rm{B}}_p(r)$ of the first kind,
\begin{align}%---------------------------------------------------------
{\rm{B}}_p(r)&\equiv\frac{{\rm{i}}^p}{2}\int_{-1}^1
\e^{-{\rm{i}}r\xi}{\rm{L}}_p(\xi)\d\xi,
\label{sBf}%-----------------------------------------------------------
\\\text{we obtain}\quad
\hat{\phi}_{\vec\jmath,k,\vec{q}}&=
\frac{1}{\pi^d}\left|\vec{\vec{h}}_k\right|
\e^{-{\rm{i}}\vec{q}\vdp\vec{a}_k}\prod_{\alpha=1}^d\sum_{p=0}^P
\check{\phi}_{\jmath^\alpha,p}{\rm{i}}^{-p}
{\rm{B}}_p(\vec{q}\vdp\vec{h}^\alpha_k).
\label{ltm}%------------------------------------------------------------
\end{align}%------------------------------------------------------------
Note that most expressions in \eqref{ltm} can be precomputed;
objects that may vary during a dynamically adaptive computation,
such as $\vec{a}_k$ or $\vec{h}^\alpha_k$, typically take
values from a sparse set, e.g., a collection of powers of 2.
%-----------------------------------------------------------------------
The computation of \eqref{Fsc} now incurs
no additional error beyond that of \eqref{ppr}.
%-----------------------------------------------------------------------
Also note, to generalize to the case $P=P_k^\alpha$ is straightforward.

\section{Accuracy of transform for 1D \& 2D test cases}
\label{rslts}

\begin{figure}\begin{center}%-------------------------------------------
\includegraphics[height=.25\textheight,width=.5\textwidth]{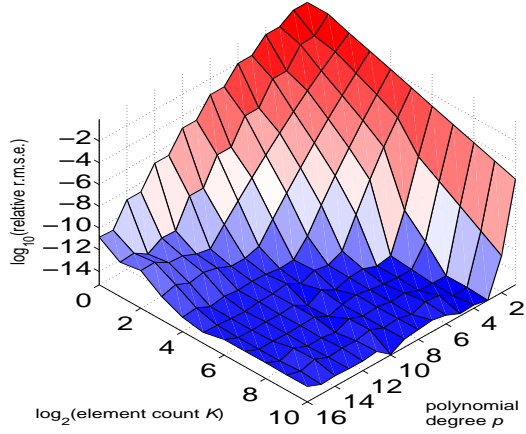}%
\caption{Surface plot (blue low to red high) of
log$_{10}$ relative r.m.s.\ error in \eqref{Fsc}
for $u^{\rm{ex}}(x)=\sin x$, vs $\log_2K$ and $P$.}
\label{f:Gqn2Fwnt}\end{center}\end{figure}

Equation \eqref{ltm} was implemented in MatLab$^{\circledR}$
and tested using known results for \eqref{Fsc}.
%-----------------------------------------------------------------------
The most immediate test follows from \eqref{sBf}, namely
$\hat{u}^{\rm{ex}}_q=\widehat{\rm{L}_p(\cdot/\pi)}_q
={\rm{i}}^{-p}{\rm{B}}_p(\pi q)$.
%-----------------------------------------------------------------------
In this case \eqref{Fsc} was found to reproduce \eqref{sBf}
to 12-16 digits for $K=1$, $P\leq18$, implying similar
performance for any polynomial $u(\vec{x})$ in this range.
%-----------------------------------------------------------------------
The next test was to put $u^{\rm{ex}}(x)=\sin x$, or
$\hat{u}^{\rm{ex}}_q=(\delta_{q,1}-\delta_{q,-1})/2{\rm{i}}$.
%-----------------------------------------------------------------------
Since this is not a polynomial we should expect at best
to see algebraic convergence w.r.t.\ $K$ in a uniform
meshing $a_k=(k-1)h_k-\pi$, $h_k=2\pi/K$ and exponential
convergence w.r.t.\ $P$, as verified in Fig.\ \ref{f:Gqn2Fwnt}.
%-----------------------------------------------------------------------
Note there is no need to test $u^{\rm{ex}}(x)=\sin rx$ for $r>1$
because of scaling.

%-----------------------------

\begin{figure}\begin{center}%-------------------------------------------
\includegraphics[width=.5\textwidth]{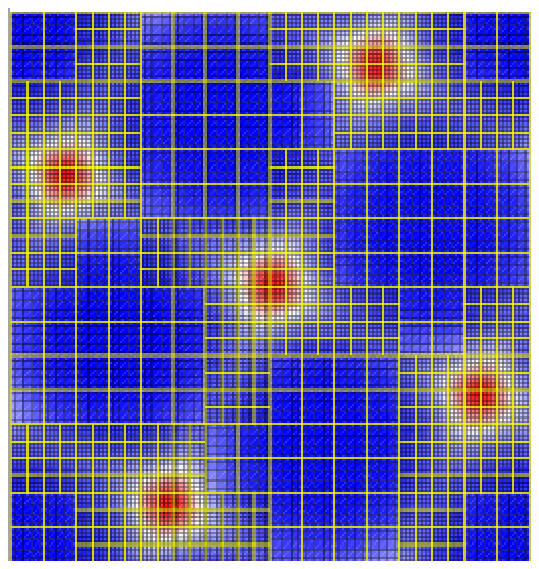}%---------------------
\includegraphics[height=.5\textwidth,width=.5\textwidth]{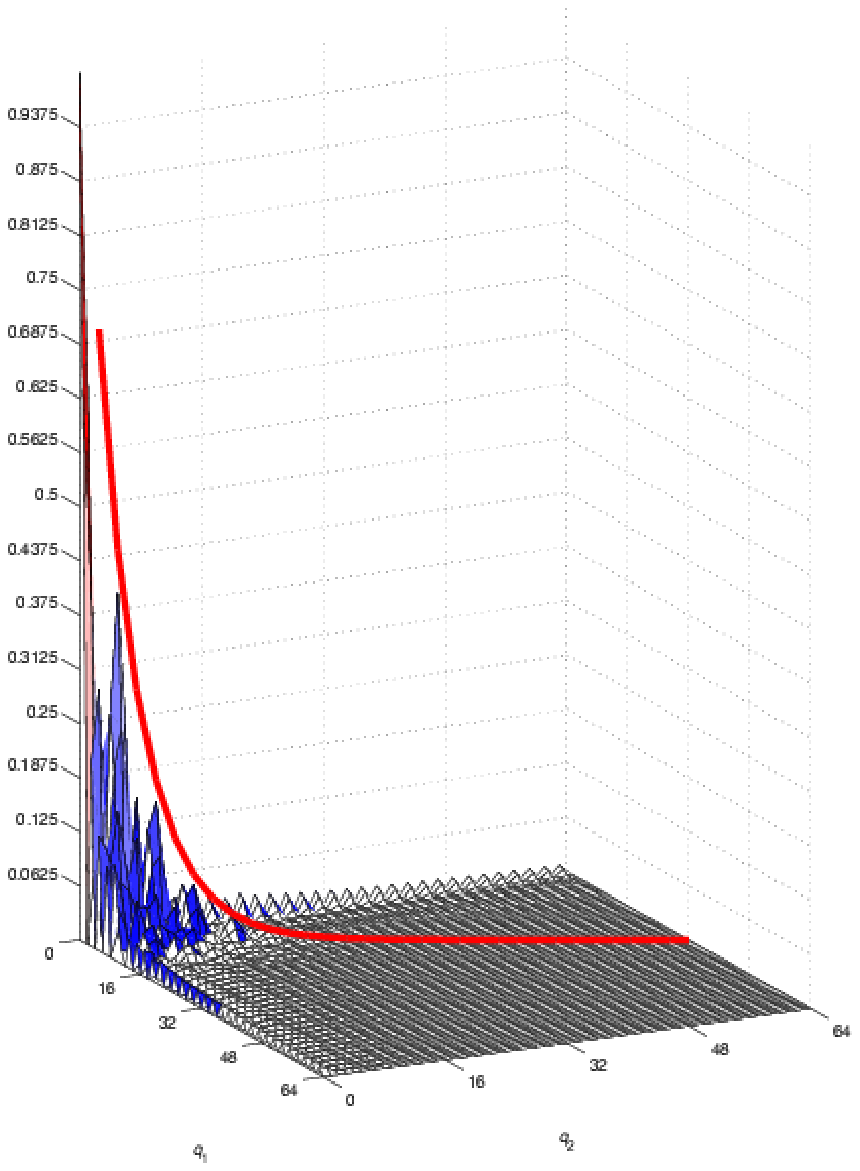}%-
\caption{Left, $u$ \eqref{e:GqnsiLor} over the spatial $\vec{x}$
domain, increasing from blue to red; yellow lines indicate element
boundaries, black lines show nodes $\vec{x}_{\vec\jmath,k}$ with $P=5$.
%-----------------------------------------------------------------------
Right, surface plot of $|\hat{u}_{\vec{q}}|$
from \eqref{Fsc} vs $q^1$ and $q^2$.}
\label{f:GqnsiLor}\end{center}\end{figure}

\begin{figure}\begin{center}%-------------------------------------------
\includegraphics[width=.5\textwidth]{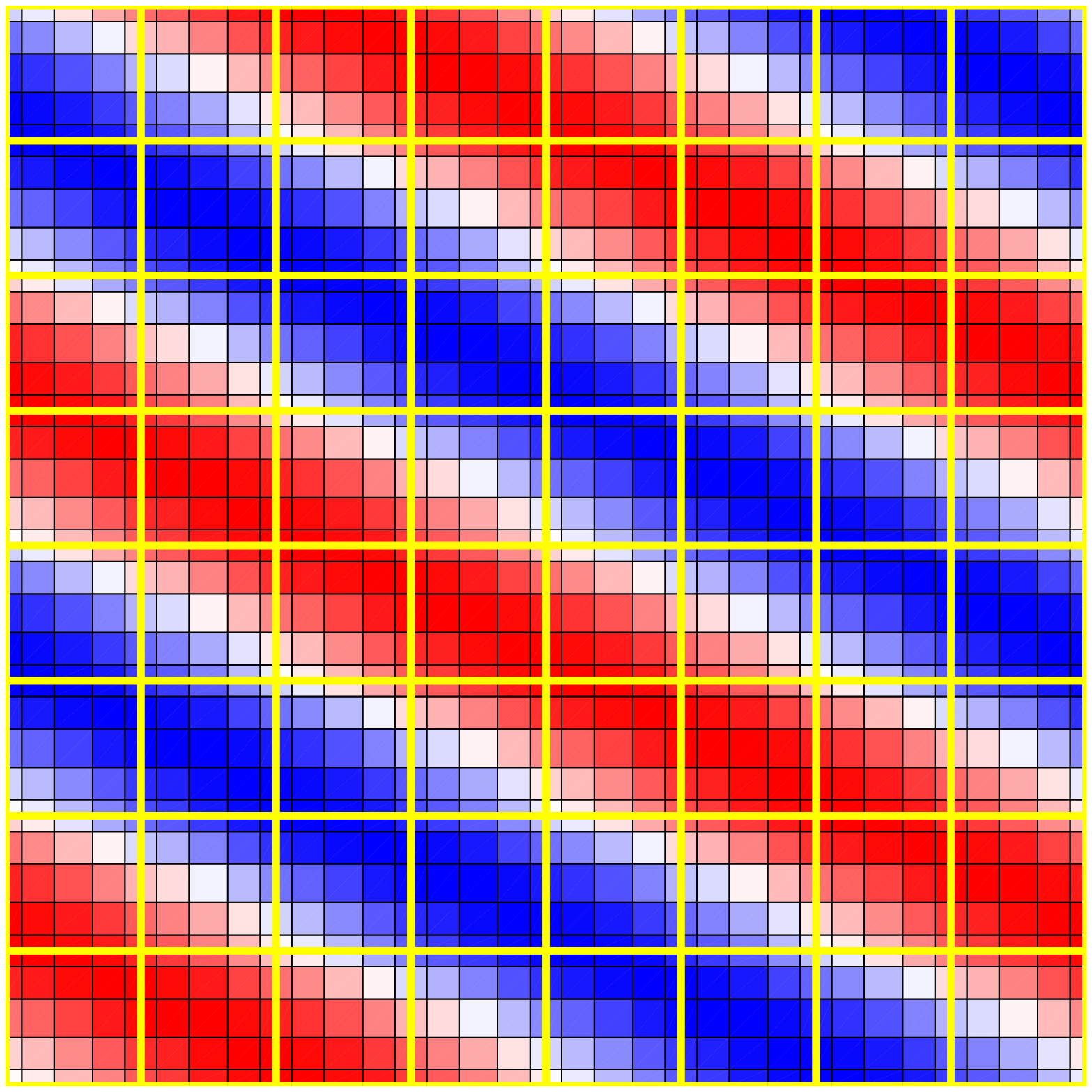}%---------------------
\includegraphics[height=.5\textwidth,width=.5\textwidth]{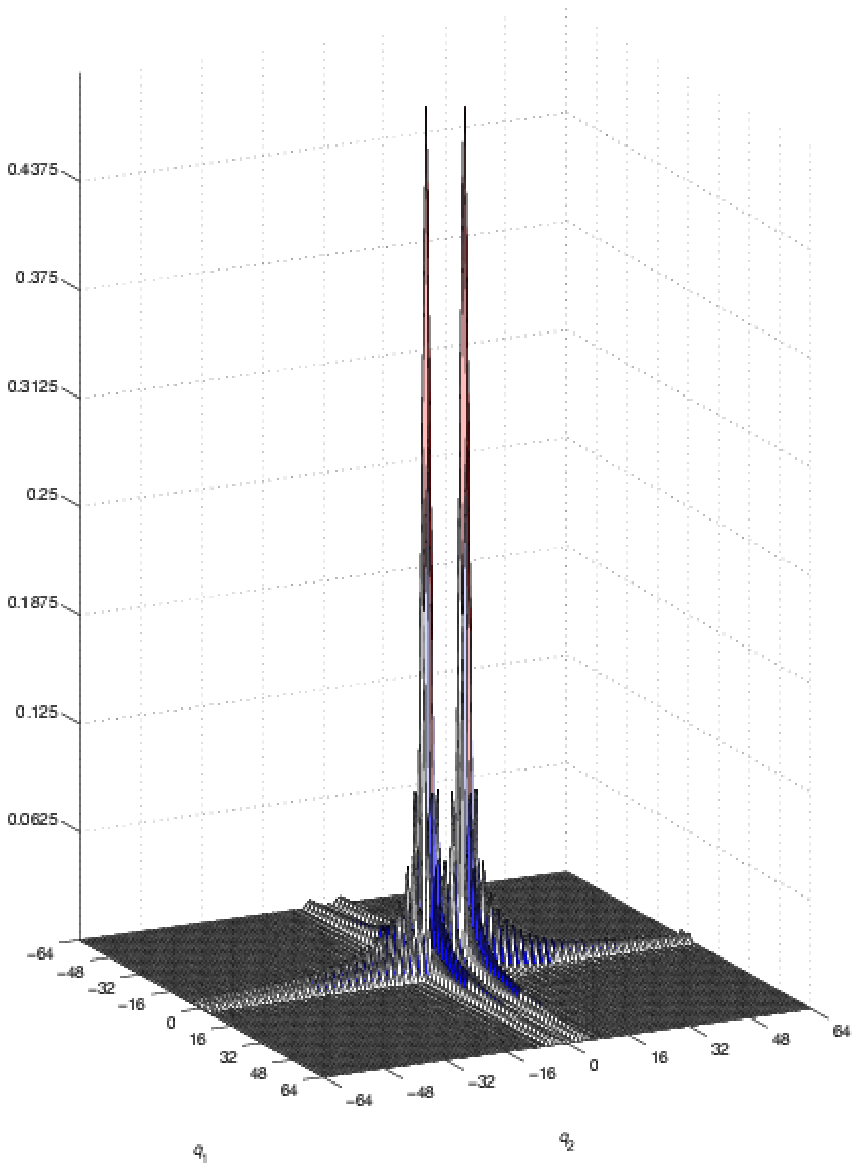}%-
\caption{As in Fig.\ \ref{f:GqnsiLor} but for the $t=0$
state given by \eqref{e:GqnBurg0}, in $K=2^6$ elements.}
\label{f:GqnBurg0}\end{center}\end{figure}

\begin{figure}\begin{center}%-------------------------------------------
\includegraphics[width=.5\textwidth]{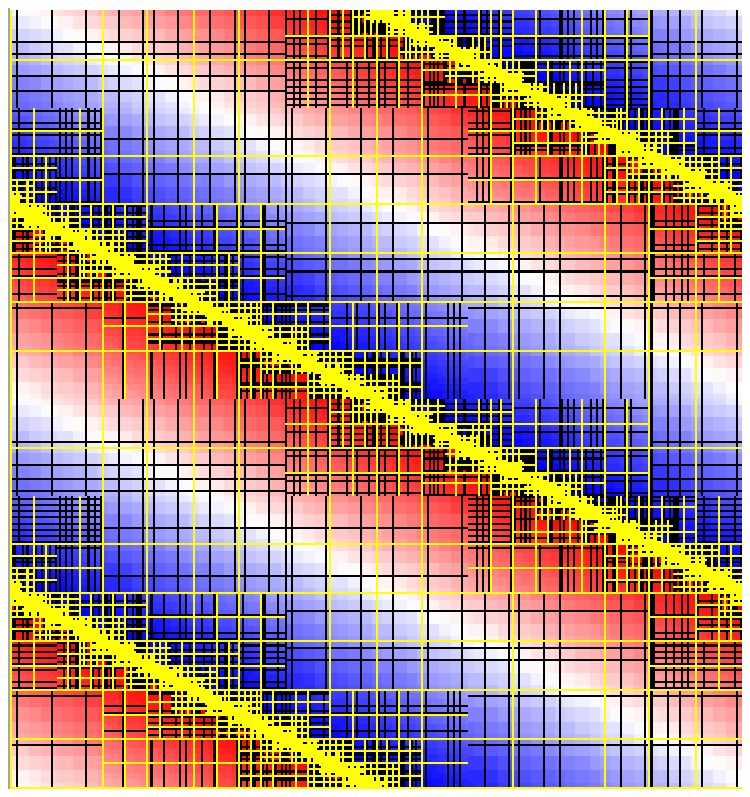}%---------------------
\includegraphics[height=.5\textwidth,width=.5\textwidth]{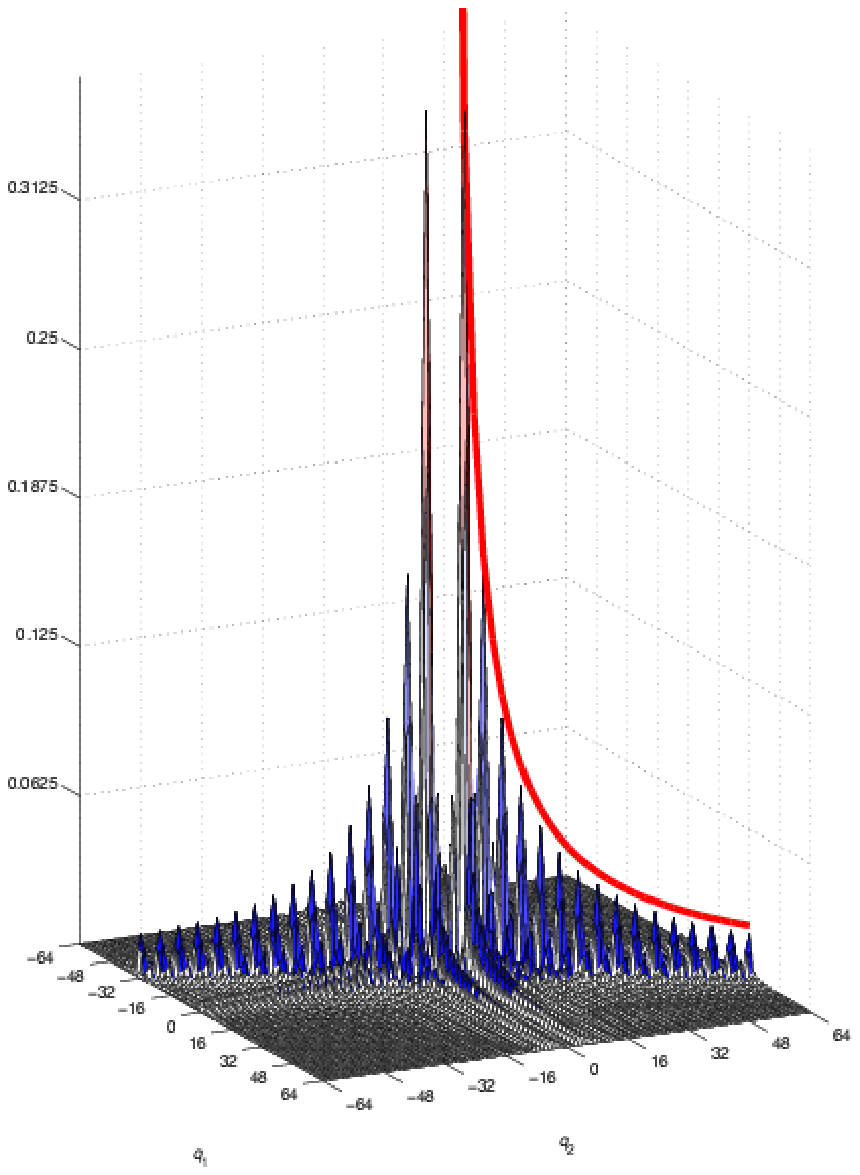}%-
\caption{As in Fig.\ \ref{f:GqnBurg0} but for $t=1.6037/5\pi$.}
\label{f:GqnBurg1}\end{center}\end{figure}

We conclude by examining three 2D tests with adaptive meshing
in the fashion of \citep{FR2005}, using MatLab$^{\circledR}$.
%-----------------------------------------------------------------------
Fig.\ \ref{f:GqnsiLor} confirms
\eqref{Fsc} in the case \citep[(19)]{fbc2005}
\begin{equation}%-------------------------------------------------------
u^{\rm{ex}}(\vec{x})\equiv\sum_{\vec{q}\in\Zset^2}
\e^{b^1|q^1|+b^2|q^2|+{\rm{i}}\vec{q}\vdp\vec{\vec{l}}\vdp\vec{x}},
\label{e:GqnsiLor}
\end{equation}%---------------------------------------------------------
where $b^\alpha=-\frac{2}{5}$ and
$\vec{\vec{l}}\doteq\left(\begin{smallmatrix}\hphantom{-}l^1&l^2\\
-l^2&l^1\end{smallmatrix}\right)=\left(\begin{smallmatrix}\hphantom{-}1&2\\
-2&1\end{smallmatrix}\right)$ is
a biperiodicity-preserving ``rotation''.
%-----------------------------------------------------------------------
As expected, the red curve (connecting the $|\hat{u}_{\vec{q}}|$
peaks) shows a power-law decay in $\vec{q}$-space.
%-----------------------------------------------------------------------
Note, in this plot and those below the $\vec{\vec{l}}$-operation
helps instigate mesh adaption but has the consequence
of leaving $\vec{q}$ undersampled in $\Zset^2$.
%-----------------------------------------------------------------------
In Fig.\ \ref{f:GqnBurg0} is shown
an initial condition \citep[(22)]{fbc2005}
\begin{equation}%-------------------------------------------------------
\vec{u}^{\rm{ex}}(0,\vec{x}):=-\vec{l}\sin\vec{l}\vdp\vec{x}
\label{e:GqnBurg0}
\end{equation}%---------------------------------------------------------
for the 2D Burgers eq.
%-----------------------------------------------------------------------
As expected, $\hat{u}_{\vec{q}}$
almost vanishes for $\vec{q}\neq\pm\vec{l}$.
%-----------------------------------------------------------------------
Finally, at time $t=1.6037/\pi|\vec{l}|^2$ the
analytic solution generalizing \citep[(2.5)]{BDHLOPO}
to 2D is shown in Fig.\ \ref{f:GqnBurg1}.
%-----------------------------------------------------------------------
As expected for the \emph{nearly} $\mathbb{C}^0$-discontinuous
fronts $\perp\vec{l}$ seen at left, $|\hat{u}^1_{\vec{q}}|$ decays
slightly faster than $\mathcal{O}(|{\vec{q}}|^{-1})$ but \emph{only
for wavevectors} $\vec{q}\|\vec{l}$ (red curve).

%-----------------------------

%-----------------------------

%-----------------------------

% Bibliographic references with the natbib package:
% Parenthetical: \citep{Bai92} produces (Bailyn 1992).
% Textual: \citet{Bai95} produces Bailyn et al. (1995).
% An affix and part of a reference:
%   \citep[e.g.][Ch. 2]{Bar76}
%   produces (e.g. Barnes et al. 1976, Ch. 2).

\end{document}